\documentclass{amsart}
\usepackage{amsfonts, amsmath, latexsym, epsfig}
\usepackage{amssymb}
\usepackage{epsf}
\usepackage{url}
\def\ml{l\kern-.15em'\kern-.085em}

\newcommand{\ZZ}{\ensuremath{\mathbb{Z}}}

\newtheorem{theorem}{Theorem}

\newtheorem{lemma}{Lemma}

\newtheorem{conjecture}{Conjecture}

\newtheorem{definition}{Definition}
\def\QuotS#1#2{\leavevmode\kern-.0em\raise.2ex\hbox{$#1$}\kern-.1em/\kern-.1em\lower.25ex\hbox{$#2$}}

\begin{document}

\author{Mathieu Dutour Sikiri\'c}
\address{Mathieu Dutour Sikiri\'c, Rudjer Boskovi\'c Institute, Bijenicka 54, 10000 Zagreb, Croatia}
\email{mdsikir@irb.hr}

\title{A variation on the Rubik's cube}

\date{}

\begin{abstract}
The Rubik's cube is a famous puzzle in which faces can be moved and the corresponding movement operations define a group.
We consider here a generalization to any $3$-valent map. We prove an upper bound on the size of the corresponding
group which we conjecture to be tight.
\end{abstract}

\maketitle

\section{Introduction}
The Rubik's cube is a $3$-dimensional toy in which each face of the cube is movable. There has been extensive study
of its mathematics (see \cite{JoynerAdventureGroup}).
On the play side the Rubik's Cube has led to the creation of many different variants (Megaminx, Pyraminx, Tuttminx,
Skewb diamond, etc.)
The common feature of those variants is that they are all physical and built as toys.

Our idea is to extend the original rubik's cube to any $3$-valent map $M$ on any surface.
To any face of the map we associate one transformation.
The full group of such transformation is named $Rubik(M)$ and we study its size and its natural normal subgroups.

A well studied extension \cite{SolvingRubikOptimally,OnTheNNNrubikCube,RubikRevenge1,RubikRevenge2} of the
Rubik's cube is the $n\times n\times n$-cube where the $3$-lanes of the cubes are extended to $n$. There has also been
interest \cite{ZeroKnowledgeRubik} in cryptographic applications of Rubik's Cube. Thus the large class of groups
that we build could be of wide interest in computer science.

In Section \ref{Construction_Rubik_map}, we construct the Rubik's cube transformation of the map.
In Section \ref{GAP_presentation}, we explain how to use existing computer algebra systems such as {\tt GAP} in
order to work with the Rubik's groups considered.
In Section \ref{Subgroup_RubikCubeGroup} we define several subgroup of the Rubik group and prove a bound on the
size of the Rubik group in the oriented case.
In Section \ref{SEC_Extensions} several possible extensions and further works are mentioned.

\section{The construction of the Rubik group from a $3$-valent map $M$}\label{Construction_Rubik_map}
In this section we use $M$ to denote a $3$-valent map.
By $V(M)$, $E(M)$, respectively $Face(M)$, we define the set of vertices, edges, respectively faces of $M$.

\begin{figure}
\begin{center}
\begin{minipage}[b]{4.3cm}
\resizebox{33mm}{!}{\rotatebox{0}{\includegraphics{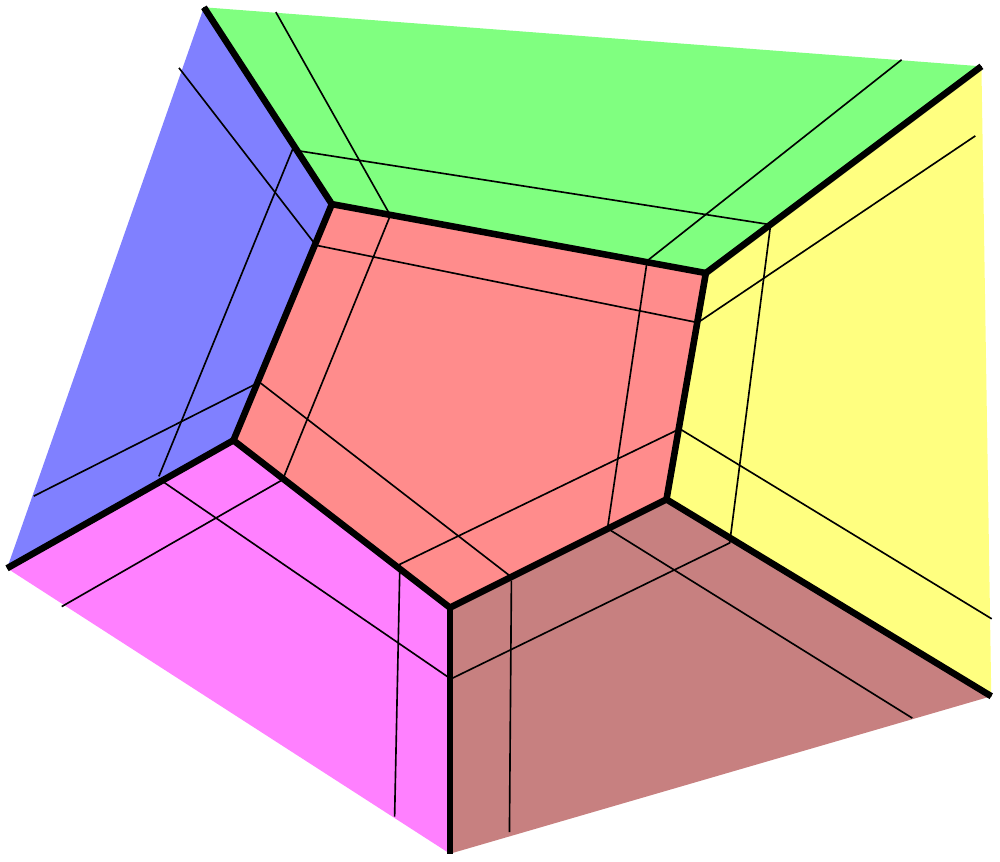}}}\par
\end{minipage}
\begin{minipage}[b]{4.3cm}
\resizebox{33mm}{!}{\rotatebox{0}{\includegraphics{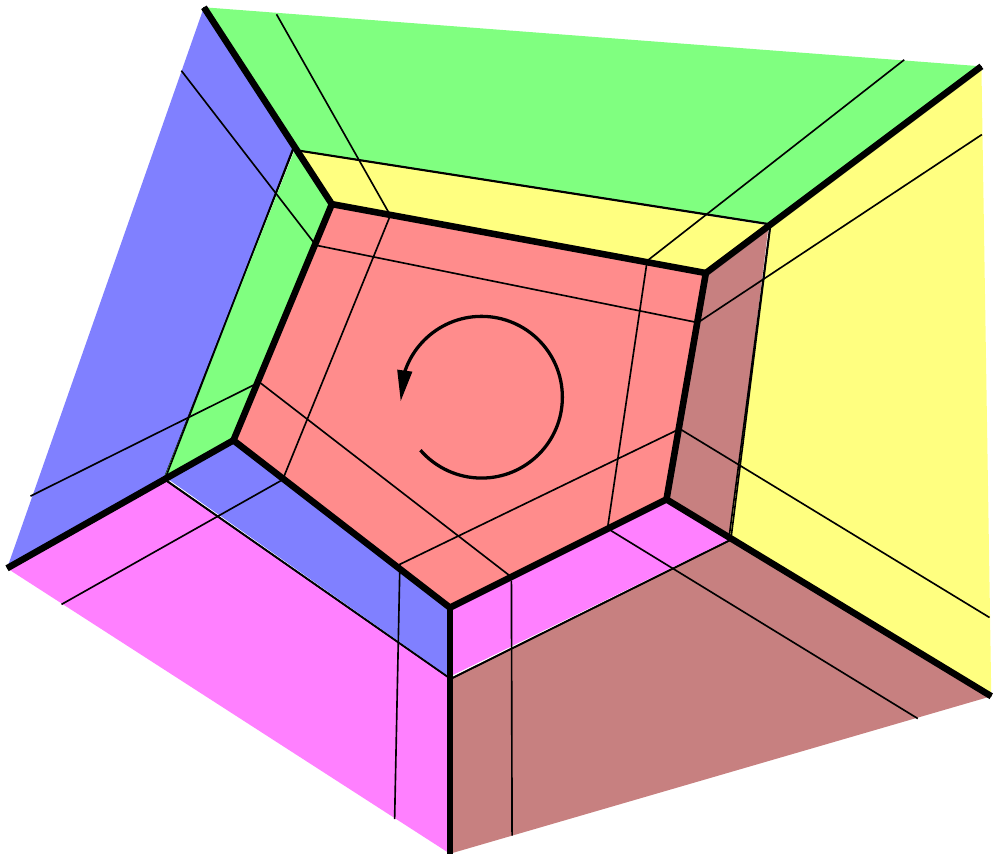}}}\par
\end{minipage}
\end{center}
\caption{The action of the side movement}
\label{SingleMovement}
\end{figure}

For any face $F$ of $M$ we define a {\em side movement} $sm(M,F)$ to be a movement in one direction of a face $F$.
This is illustrated in Figure \ref{SingleMovement}. If the face contains $p$ edges then the element $sm(M, F)$
has order $p$.

For a set $S$ of elements of a group $G$ we define $Group(S)$ to be the group generated by
$S$.
Using this we can define the Rubik's group:
\begin{equation*}
Rubik(M) = Group(\left\{sm(M, F) \mbox{~for~} F\in Face(M)\right\})
\end{equation*}

\begin{figure}
\begin{center}
\begin{minipage}[b]{5.5cm}
\resizebox{33mm}{!}{\rotatebox{0}{\includegraphics{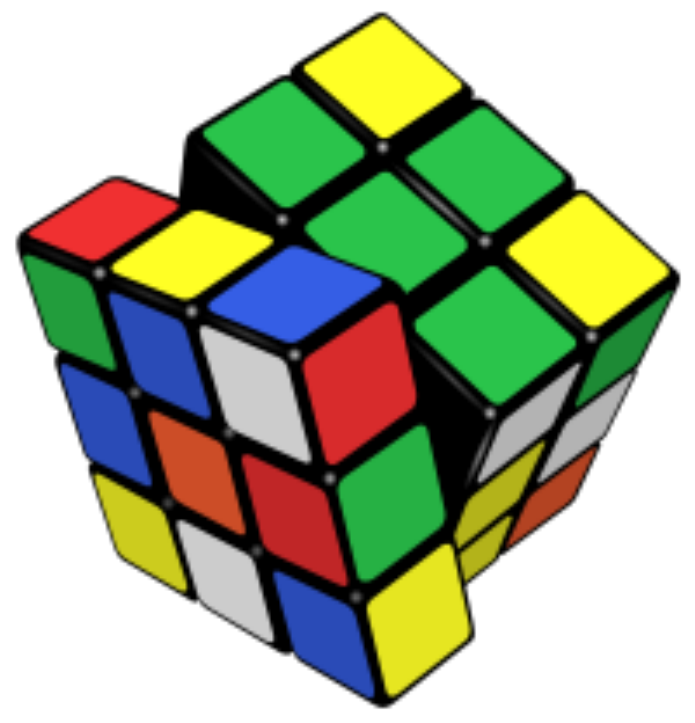}}}\par
$Rubik(Cube)$: Rubik's Cube.
\end{minipage}
\begin{minipage}[b]{5.5cm}
\resizebox{33mm}{!}{\rotatebox{0}{\includegraphics{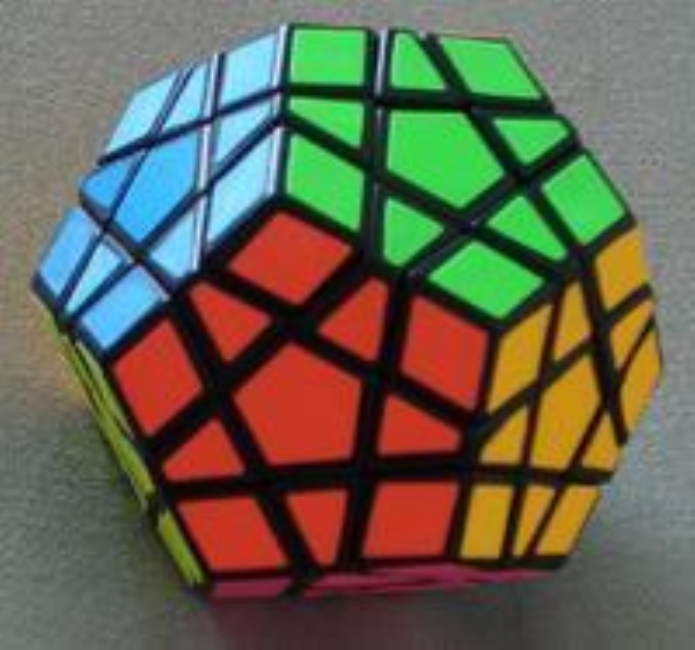}}}\par
$Rubik(Dodecahedron)$: Megaminx
\end{minipage}
\end{center}
\caption{The Rubik construction for Cube and Dodecahedron}
\label{RubicPlatonic}
\end{figure}

For the Cube and Dodecahedron this is known as Rubik's cube and Megaminx, see Figure \ref{RubicPlatonic}
where the toy can be represented physically.

A {\em corner} of a map $M$ is a pair $(F, v)$ with $v$ a vertex contained in a face $F$.
A {\em side edge} of a map $M$ is a pair $(F, e)$ with $e$ an edge contained in a face $F$.
A $3$-regular map $M$ with $v$ vertices has exactly $e = 3v / 2$ edges.
It will have $3v$ corners and $2e = 3v$ side edges.
By $Corner(M)$, we define the set of corners of $M$.

The face movement $sm(M, F)$ acts on the set of corners and the set of side edges.
This corresponds directly to the faces of the Rubik's cube.

\section{The computer Algebra side of things}\label{GAP_presentation}

We  present the group $Rubik(M)$ as a permutation group on the corners and side edges.
This way of presenting the group follows \cite{ZassenhausRubik}
\footnote{See also an example in \url{https://www.gap-system.org/Doc/Examples/rubik.html}}
where the standard Rubik's cube was considered.

\begin{figure}
\begin{center}
\begin{minipage}[b]{6.3cm}
\resizebox{43mm}{!}{\rotatebox{0}{\includegraphics{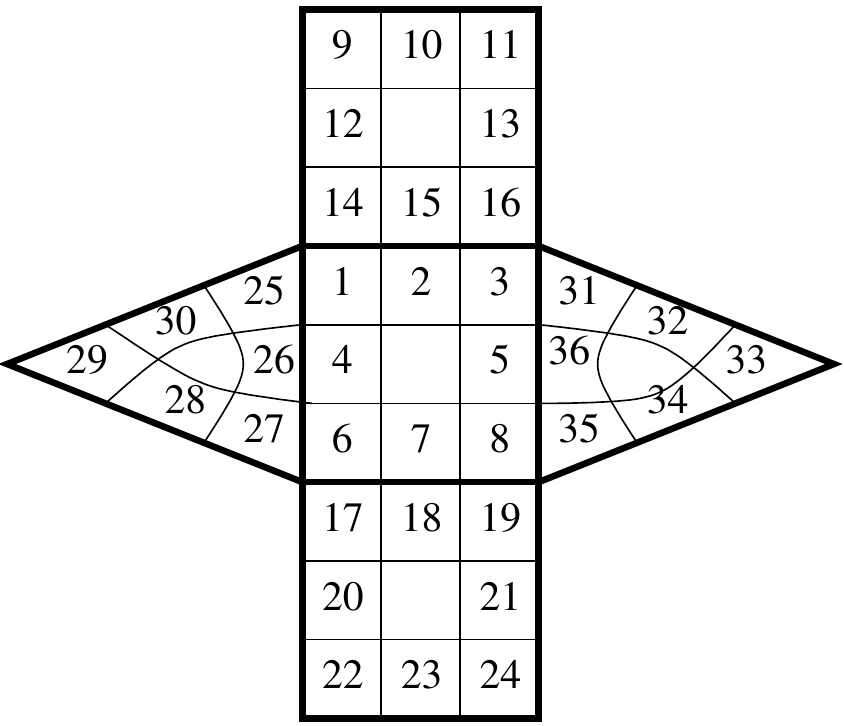}}}\par
\end{minipage}
\end{center}
\caption{Numbering of corners and side-edge of $Prism_3$}
\label{Numbering_Prism3}
\end{figure}

For the example of $Prism_3$ plane graph, we can number the corners and side-edges according to
Figure \ref{Numbering_Prism3}. The side-movements generators of $Rubik(Prism_3)$ can then be
expressed as the following permutations:

{\scriptsize
\begin{verbatim}
gap> RubikPrism3:=Group([(31,35,33)(36,34,32)(3,19,11)(5,21,13)(8,24,16),
  (28,27,25)(30,28,26)(6,14,22)(4,12,20)(1,9,17),
  (1,6,8,3)(4,7,5,2)(17,35,16,25)(19,31,14,27)(18,36,15,26),
  (19,17,22,24)(18,20,23,21)(7,28,10,34)(8,27,9,33)(6,29,11,35),
  (9,14,16,11)(12,15,13,10)(1,31,24,29)(2,32,23,30)(3,33,22,25)]);
\end{verbatim}
}

The advantage of this approach is that very efficient implementations of
permutation group algorithms ({\tt GAP}) which can address in an efficient
way many questions one may have about a permutation group.

One famous problem in Rubik's cube group theory is to find the minimal number of
moves to go from one configuration to another, in other words to find the diameter
of the Cayley graph of the group with the side moves as generators
(see \cite{DiameterRubikCubeGroup}).
The {\tt GAP} computer algebra software does not provide a solution to that problem,
but it allows to express and element of the group in terms of the generators,
that is find an expression, possibly not minimal.

More specific functionality for the computation is available in the GAP package
\cite{plangraph}.

\section{Subgroups of $Rubik(M)$}\label{Subgroup_RubikCubeGroup}

The group $Rubik(M)$ acts not only of the corners and side edges, but also on the
vertices and edges. To represent those actions we use subscript such as
$sm(M,F)_{edge}$ and denote the corresponding group $Rubik(M)_{edge}$.

\begin{lemma}\label{SignatureLemma}
  Given a map $M$ and a face $F$ we have:

  (i) The signature of $sm(M,F)_{side\, edge}$ is $1$.

  (ii) The signature of $sm(M,F)_{vertex}$ is equal to the signature of $sm(M,F)_{edge}$.

\end{lemma}
\proof (i) If $F$ is a face with $p$ sides then $F$ contains $p$ side edges $(F,e)$,
each side edges is adjacent to another side edge $(F_e, e)$.
The side generator $sm(M,F)$ has one cycle formed by the side-edges $(F,e)$ and another cycle
formed by the side edges $(F_e, e)$. Each of those cycles has the same length and so the signature
of $sm(M,F)$ acting on the side-edges is $1$.

(ii) If $F$ is a face with $p$ sides then it contains $p$ edges and $p$ vertices. The action on
each of them is a cycle of length $p$. So their signature is the same. \qed

We now consider what additional properties can be obtained if the map $M$ is oriented.

\begin{definition}
Assume $M$ is an oriented map. Every vertex $v$ of $M$ is contained in $3$
corners $C_v=\{c_1, c_2, c_3\}$. Since $M$ is orientable, we can assume that $C_v$
is oriented in the direct way.

We denote by $OrMap(M)$ the set of permutation of the corners that preserve the
triples and their orientation.

Given two corners $c$, $c'$ we call $sh(c,c')$ the number $k$ of rotations $r$ by $2\pi /3$
such that $r^k(c) = c'$.
\end{definition}
The orientability of the map will allow us to restrict the possibilities of mapping the
corners. We need the $OrMap(M)$ set of map in order to define some properties below
for the Rubik's cube group. The shift is more complicated. If a vertex $v$ has $C_v = \{c_1, c_2, c_3\}$
in direct orientation then we have e.g. $sh(c_i, c_i) = 0$, $sh(c_1, c_2) = 1$
and $sh(c_1, c_3) = 2$.

\begin{theorem}\label{ThreeOrientability}
Let $M$ be an oriented map. The following holds:

(i) For $c$, $c'$, $c''$ corners we have $sh(c, c'') = sh(c,c') + sh(c',c'')\pmod 3$

(ii) For a function $\phi:V(M) \rightarrow Corn(M)$ such that $\phi(v)$ is a corner of $V$
and a function $f\in OrMap(M)$ the expression
\begin{equation*}
\sum_{v\in V(M)} sh(f(\phi(v)), \phi(f(v)))
\end{equation*}
is independent of $\phi$ and denoted $sh(f)$.

(iii) If $f, g \in OrMap(M)$ then we have
\begin{equation*}
sh(f\circ g) = sh(f) + sh(g) \pmod 3
\end{equation*}

(iv) $sh$ is a non-trivial function on $OrMap(M)$.

(v) For $f\in Rubik(M)$ we have $sh(f) = 0\pmod 3$.

\end{theorem}
\proof (i) The equality follows by additivity of the rotational shift.

(ii) Let $\phi$, $\phi'$ two such functions. Then we define $r(v) = sh(\phi(v), \phi'(v))$.
Define first $a(\phi) = \sum_{v\in V(M)} sh(f(\phi(v)), \phi(f(v)))$, we then have:
\begin{equation*}
\begin{array}{rcl}
a(\phi')
&=& \sum_{v\in V(M)} sh(f(\phi'(v)), f(\phi(v)) + sh(f(\phi(v)), \phi(f(v))\\
&& + sh(\phi(f(v)), \phi'(f(v)))\\
&=& a(\phi) + \sum_{v\in V(M)} sh(f(\phi'(v)), f(\phi(v)) + \sum_{v\in V(M)} sh(\phi(f(v)), \phi'(f(v)))\\
&=& a(\phi) + \sum_{v\in V(M)} sh(\phi'(v)), \phi(v)) + \sum_{v\in V(M)} sh(\phi(f(v)), \phi'(f(v)))\\
&=& a(\phi) + \sum_{v\in V(M)} sh(\phi'(v)), \phi(v)) + \sum_{v\in V(M)} sh(\phi(v)), \phi'(v))\\
&=& a(\phi) + \sum_{v\in V(M)} sh(\phi'(v)), \phi(v)) + sh(\phi(v)), \phi'(v))\\
&=& a(\phi) + \sum_{v\in V(M)} sh(\phi'(v)), \phi'(v))\\
&=& a(\phi)
\end{array}
\end{equation*}
The equality $sh(f(\phi'(v)), f(\phi(v))) = sh(\phi'(v), \phi(v))$ comes from the fact that $f$
preserves the orientation of the corners.
The equality $\sum_{v\in V(M)} sh(\phi(f(v)), \phi'(f(v))) = \sum_{v\in V(M)} sh(\phi(v), \phi'(v))$
comes from the fact that $f$ is permuting the vertices of $M$.

(iii) Let $f$ and $g$ be two such mappings. We then have
\begin{equation*}
\begin{array}{rcl}
sh(f\circ g)
&=& \sum_{v\in V(M)} sh(f(g(\phi(v))), \phi(f(g(v))))\\
&=& \sum_{v\in V(M)} sh(f(g(\phi(v))), f(\phi(g(v)))) + sh(f(\phi(g(v))), \phi(f(g(v))))\\
&=& \sum_{v\in V(M)} sh(g(\phi(v)), \phi(g(v))) + \sum_{v\in V(M)} sh(f(\phi(v)), \phi(f(v)))\\
&=& sh(g) + sh(f)
\end{array}
\end{equation*}
We have used that $f$ preserves orientation and the fact that $g$ permutes the vertices in our rewriting.

(iv) Let us select a vertex $v$ of $M$ with $C_v= \{c_1, c_2, c_3\}$. We set $f(c_i) = c_{i+1}$
for $c_i\in C_v$ and $f(c) = c$ otherwise. We find that $sh(f) = 1$.

(v) Let $F$ be a face of $M$. We select a function $\phi$ such that $phi(v) = (F,v)$ for $v$
and vertex contained in $F$ and an arbitrary choice on other vertices.
Setting $f = sm(M,F)$ we have
\begin{equation*}
sh(f) = \sum_{v\in F\cap V(M)} sh(f(\phi(v)), \phi(f(v)))
\end{equation*}
It is easy to see that $\phi(f(v)) = f(\phi(v))$ for $v\in F\cap V(M)$.
Thus we get $sh(f)$. Since the side movement generate $Rubik(M)$, by (iii) we get that for all
$f\in Rubik(M)$, $sh(f) = 0$. \qed

\begin{theorem}
We have a sequence of group homomorphisms:
\begin{equation*}
\begin{array}{c}
  Rubik(M) = Rubik(M)_{corner, side\,\,edge} \rightarrow Rubik(M)_{corner, edge}\\
  \mbox{~~~~~~~~} \rightarrow Rubik(M)_{corner} \rightarrow Rubik(M)_{vertex}
\end{array}
\end{equation*}

\end{theorem}
\proof The equality $Rubik(M) = Rubik(M)_{corner, side\,\,edge}$ comes from the definition
of the Rubik group.
If the group acts on the side-edges, then it acts on the edges and so we get the map
$Rubik(M)_{corner, side\,\,edge} \rightarrow Rubik(M)_{corner, edge}$.
Dropping the edge action gets us the mapping to $Rubik(M)_{corner}$.
Finally, a vertex is contained in $3$ faces and so in $3$ corners.
This defines a mapping from $Rubik(M)_{corner}$ to $Rubik(M)_{vertex}$. \qed

We define following subgroups from those homomorphisms:
\begin{equation*}
  \left\{\begin{array}{rcl}
  H_1(M) &=& Ker(Rubik(M)_{corner, side\,\,edge} \rightarrow Rubik(M)_{corner, edge})\\
  H_2(M) &=& Ker(Rubik(M)_{corner, edge} \rightarrow Rubik(M)_{corner})\\
  H_3(M) &=& Ker(Rubik(M)_{corner} \rightarrow Rubik(M)_{vertex})
  \end{array}\right.
\end{equation*}

We can now use the above decomposition in order to get a conjectural description of $Rubik(M)$.
\begin{theorem}
  Given an oriented map $M$ the following holds:

  (i) $H_1(M)$ is a subgroup of $\ZZ_2^{|E(M)|-1}$

  (ii) $H_2(M)$ is a subgroup of $A_{|E(M)|}$.

  (iii) $H_3(M)$ is a subgroup of $\ZZ_3^{|V(M)|-1}$

  (iv) $Rubik(M)_{vertex}$ is a subgroup of $A_{|V(M)|}$ if all faces of $M$ have odd size and a subgroup of $S_{|V(M)|}$ otherwise.

\end{theorem}
\proof $H_1(M)$ is formed by all the transformations that preserve all corners and edges
but may switch a side edge $(F_1, e)$ into another side edges $(F_2, e)$.
Since every edge is contained in two faces, this makes $H_1(M)$ a subgroup of the commutative
subgroup $\ZZ_2^{|E(M)|}$.
The group $\ZZ_2^{|E(M)|}$ contains some elements that switch
just two side-edges and are thus of signature $-1$. Thus by Lemma \ref{SignatureLemma}.(i) $H_1(M)$
is isomorphic to a strict subgroup of $\ZZ_2^{|E(M)|}$.

$H_2(M)$ is formed by the transformations that preserves the corners but will permutes the edges.
Thus we have that $H_2(M)$ is a subgroup of $S_{|E(M)|}$.
Since elements of $H_2(M)$ preserves all vertices, the signature of their action is $1$.
By Lemma \ref{SignatureLemma}.(ii) the signature of their action on edges is also $1$.
Therefore $H_2(M)$ is a subgroup of $A_{|E(M)|}$.

$H_3(M)$ is formed by all transformations that preserves the vertices but may permutes the corners.
Since the corners are oriented, the operation on each vertex $v$ may be encoded by an element $x_v$
of $\ZZ_3$. By Theorem \ref{ThreeOrientability}.(v) we have $sh(f)=0$ for each $f\in H_3(M)$.
Thus we get $\sum_{v\in V(M)} x_v = 0$. This means that $H_3(M)$ is a subgroup of $\ZZ_3^{|V(M)|-1}$.

The side movements $sm(M, F)$ act on the vertices. For a face of size $p$ the signature is $(-1)^{p-1}$.
Thus if all faces of $M$ have odd size then all side movements have signature $1$ and $Rubik(M)_{vertex}$
is a subgroup of $A_{|V(M)|}$. If there is a face of even size then it is a subgroup of $S_{|V(M)|}$. \qed

\begin{conjecture}
(``First law of cubology'') Given an oriented map $M$ with $v$ vertices and $e$ edges the following holds:

  (i) $H_1(M)$ is isomorphic to $\ZZ_2^{|E(M)|-1}$

  (ii) $H_2(M)$ is isomorphic to $A_{|E(M)|}$.

  (iii) $H_3(M)$ is isomorphic to $\ZZ_3^{|V(M)|-1}$

  (iv) $Rubik(M)_{vertex}$ is isomorphic to $A_{|V(M)|}$ if all faces of $M$ have odd size and $S_{|V(M)|}$ otherwise.

\end{conjecture}

This conjecture has been checked for many plane graphs e.g. the ones with at most $40$ vertices and faces of size $6$ or $p$ with $3\leq p\leq 5$.
This Theorem is proved for the cube case in \cite[Theorem 11.2.1]{JoynerAdventureGroup},
\cite[Section 2.4]{Bandelow} and
\cite[Theorem 1.3.24]{TechRepKTH}.

\section{Possible extensions and open problems}\label{SEC_Extensions}

\begin{figure}
\begin{center}
\begin{minipage}[b]{4.3cm}
\resizebox{33mm}{!}{\rotatebox{0}{\includegraphics{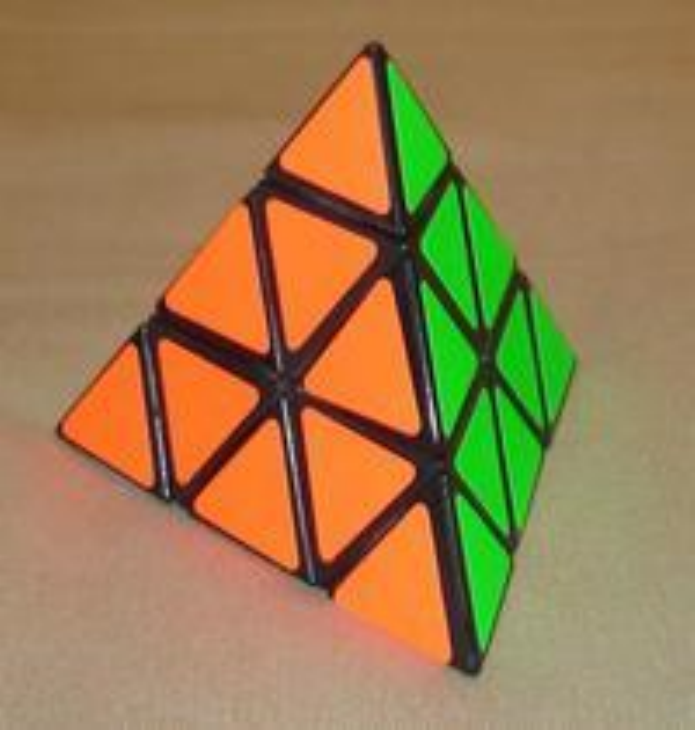}}}\par
\end{minipage}
\end{center}
\caption{The Pyraminx toy}
\label{Pyraminx}
\end{figure}

There are many possible extensions of this work. First we could consider other constructions
so as to generalize all the existing Rubik's cube variant toys to a combinatorial setting.
For example the Pyraminx (Figure \ref{Pyraminx}) does not belong to the family described here.

Proving the conjecture on the description of $Rubik(M)$ appears to be a difficult problem as
the cube proofs seem difficult to generalize.
Another research question is to consider the classification
of elements of order $2$ of $Rubik(M)$ as well as the determination of the maximal
order of the elements of $Rubik(M)$ (see \cite{TechRepKTH} for the Cube case).

Much time and energy has been spent on playing on Rubik's cube and this could be done as well
for $Rubik(M)$.
One would have to forget the physical part and accept playing on a smartphone, which could be
programmed reasonably easily. For plane graph this could be represented by a Schlegel diagram
and for toroidal maps, we could use a plane representation with a group of symmetries.

Classifying the $3$-valent graphs for which $Rubik(M)$ may be constructed physically
is also an interesting problem.

\bibliographystyle{amsplain_initials_eprint}
\bibliography{RefRubik}

\providecommand{\bysame}{\leavevmode\hbox to3em{\hrulefill}\thinspace}
\providecommand{\MR}{\relax\ifhmode\unskip\space\fi MR }
\providecommand{\MRhref}[2]{%
  \href{http://www.ams.org/mathscinet-getitem?mr=#1}{#2}
}
\providecommand{\href}[2]{#2}
\begin{thebibliography}{10}

\bibitem{Bandelow}
C.~Bandelow, \emph{Inside {R}ubik's cube and beyond}, Birkh\"{a}user, Boston,
  Mass., 1982, Translated from the German by Jeannette Zehnder [Jeannette
  Zehnder-Reitinger] and Lucy Moser.

\bibitem{TechRepKTH}
O.~Bergvall, E.~Hynning, M.~Hedberg, J.~Mickelin, and P.~Masawe, \emph{On
  {R}ubik's cube}, Report, KTH Royal Institute of Technology, 2010, URL:
  \url{https://people.kth.se/~boij/kandexjobbVT11/Material/rubikscube.pdf}.

\bibitem{RubikRevenge1}
S.~Bonzio, A.~Loi, and L.~Peruzzi, \emph{The first law of cubology for the
  {R}ubik's revenge}, Math. Slovaca \textbf{67} (2017), no.~3, 561--572.

\bibitem{OnTheNNNrubikCube}
S.~Bonzio, A.~Loi, and L.~Peruzzi, \emph{On the {$n\times n\times n$} {R}ubik's
  cube}, Math. Slovaca \textbf{68} (2018), no.~5, 957--974.

\bibitem{SolvingRubikOptimally}
E.~D. Demaine, S.~Eisenstat, and M.~Rudoy, \emph{Solving the {R}ubik's cube
  optimally is {NP}-complete}, 35th {S}ymposium on {T}heoretical {A}spects of
  {C}omputer {S}cience, LIPIcs. Leibniz Int. Proc. Inform., vol.~96, Schloss
  Dagstuhl. Leibniz-Zent. Inform., Wadern, 2018, pp.~Art. No. 24, 13.

\bibitem{plangraph}
M.~Dutour~Sikiri\'c, \emph{Plangraph},
  \url{http://mathieudutour.altervista.org/PlanGraph/}.

\bibitem{JoynerAdventureGroup}
D.~Joyner, \emph{Adventures in group theory}, second ed., Johns Hopkins
  University Press, Baltimore, MD, 2008, Rubik's cube, Merlin's machine, and
  other mathematical toys.

\bibitem{RubikRevenge2}
M.~E. Larsen, \emph{Rubik's revenge: the group theoretical solution}, Amer.
  Math. Monthly \textbf{92} (1985), no.~6, 381--390.

\bibitem{DiameterRubikCubeGroup}
T.~Rokicki, H.~Kociemba, M.~Davidson, and J.~Dethridge, \emph{The diameter of
  the {R}ubik's cube group is twenty}, SIAM J. Discrete Math. \textbf{27}
  (2013), no.~2, 1082--1105.

\bibitem{ZeroKnowledgeRubik}
E.~Volte, J.~Patarin, and V.~Nachef, \emph{Zero knowledge with {R}ubik's cubes
  and non-abelian groups}, Cryptology and network security, Lecture Notes in
  Comput. Sci., vol. 8257, Springer, Cham, 2013, pp.~74--91.

\bibitem{ZassenhausRubik}
H.~Zassenhaus, \emph{Rubik's cube: a toy, a {G}alois tool, group theory for
  everybody}, Phys. A \textbf{114} (1982), no.~1-3, 629--637.

\end{thebibliography}

\end{document}